\documentclass{amsart}
\usepackage{amsmath, amssymb, amsthm, graphicx}

\theoremstyle{plain}
\newtheorem{prop}{Proposition}[section]

\newtheorem{lemm}[prop]{Lemma}
\newtheorem{ques}[prop]{Question}

\theoremstyle{definition}
\newtheorem*{defi}{Definition}

\newtheorem{exam}[prop]{Example}
\newtheorem{rema}[prop]{Remark}

\numberwithin{equation}{section}

\def\Reff#1; #2; #3; #4; #5; #6; #7\par{%
\bibitem{#1} #2, {\it #3}, #4 {\bf #5} (#6) #7}

\def\Ref#1; #2; #3; #4\par{%
\bibitem{#1} #2, {\it #3}, #4}

\def\AA{A}
\def\AAp{\AA^+}
\def\at{\s_1}
\def\att{\s_2}

\let\b=\beta

\def\cc{s}
\def\ccc{t}
\def\cl#1{\overline{\vrule width 0pt height 6pt #1}}

\let\D=\Delta
\def\diag{D}		

\def\divel{\preccurlyeq}
\def\dr{\mathord{\backslash}}

\def\e{\varepsilon}

\def\ff{\phi}

\let\g=\gamma
\def\ge{\geqslant}
\def\gen{\sigma}
\def\genn{\tau}
\def\gennn{\rho}
\def\GG{G}
\def\GGp{\GG^+}

\def\ie{{\it i.e.}}
\def\ii{^{-1}}

\def\le{\leqslant}

\def\multer{\succcurlyeq}

\def\n(#1){\Vert #1 \Vert}

\let\pp=\dots
\def\Prod{\mbox{\rm prod}}

\def\resp{{\it resp.\ }}

\def\s{\sigma}
\def\sp{s}		
\def\spp{t}		
\def\ss#1{\sigma_{#1}^{\phantom 1}}  
\def\sss#1{\sigma_{#1}^{-1}}  
\def\SS{S}

\def\val{\nu}

\def\xx{x}

\def\yy{y}

\def\zz{z}
\def\ZZ{\mathbb Z}

\begin{document}

\author{Patrick DEHORNOY}
\address{Laboratoire de Math\'ematiques Nicolas
Oresme, CNRS UMR 6139, Universit\'e de Caen,
14032 Caen, France}
\email{dehornoy@math.unicaen.fr}
\urladdr{//www.math.unicaen.fr/\!\hbox{$\sim$}dehornoy}

\title{Disks in trivial braid diagrams}



\begin{abstract}
We show that every trivial $3$-strand braid diagram
contains a disk, defined as a ribbon ending in
opposed crossings. Under a convenient algebraic
form, the result extends to every Artin--Tits
group of dihedral type, but it fails to
extend to braids with $4$~strands and more. The
proof uses a partition of the Cayley graph and a
continuity argument.
\end{abstract}

\maketitle

\section{Introduction}\label{S:Intro}

Let us say that a braid diagram is trivial if it represents
the unit braid, \ie, if it is isotopic to an unbraided
diagram. Consider the following simple trivial
diagrams:
$$\includegraphics*[scale=0.7]{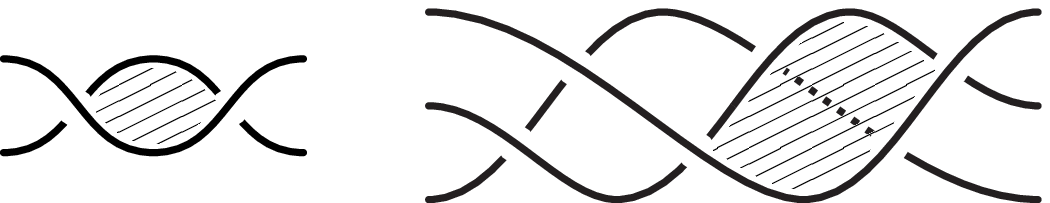}$$
We see that these digrams contain a {\it disk},
defined as an embedded ribbon ending in crossings
with opposite orientations (the striped areas). Below is
another trivial braid diagram containing a disk:
here the shape is more complicated, but we still
have the property that the third strand does not
pierce the disk.
$$\includegraphics*[scale=0.7]{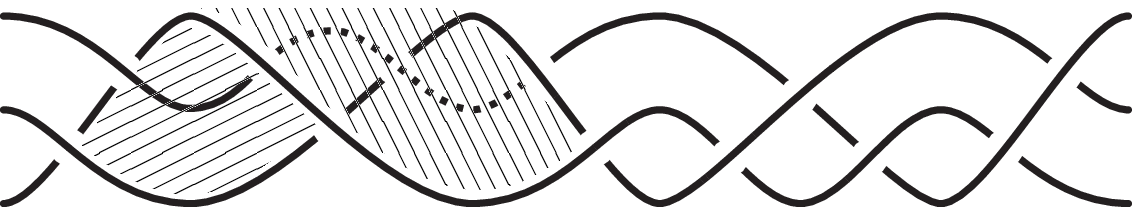}$$
Finally, let us display a more intricate example
involving a disk: here the third strand pierces the
ribbon, but it does it so as to make a topologically
trivial handle through the disk, so, up to an isotopy, we
still have an unpierced disk.
$$\includegraphics*[scale=0.65]{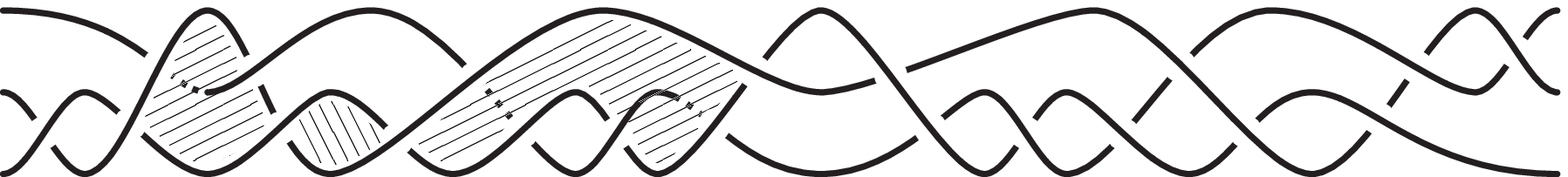}$$

A few tries should convince the reader that most trivial
braid diagrams seem to contain at least one disk in
the sense above---a precise definition will be given
below---and make the following question natural:

\begin{ques}\label{Q:1}
 Does every trivial braid diagram (with at least one
 crossing) contain a disk?
\end{ques}

Our aim is to anwser the question by proving

\begin{prop}\label{P:Pair}
The answer to Question~\ref{Q:1} is positive in the
case of $3$-strand braids, \ie, every trivial $3$-strand
braid diagram with at least one crossing contains a
disk. It is negative in the case of $4$~strands and
more.
\end{prop}

As for the negative part, it is sufficient to
exhibit a counter-example, what will be done at
the end of Section~\ref{S:Disk} (see
Figure~\ref{F:Fig2}).

As for the positive part, the argument consists in going
to the Cayley graph of the braid group and using a
continuity result, which itself relies on the properties
of division in the braid monoid~$B_n^+$. The argument
works in every Artin--Tits group of spherical type,
and we actually prove the counterpart of (the positive
part of) Proposition~\ref{P:Pair} in all Artin--Tits groups
of type~$I_2(m)$.

One should keep in mind that we are interested in
braid diagrams, not in braids: up to an isotopy, all braid
diagrams we consider can be unbraided. What makes
the question nontrivial is that isotopy may change
the possible disks of a braid diagram completely,
so that it is hopeless to trace the disks along an
isotopy. For instance, the reader can check that
applying one type~III Reidemeister move in the
braid diagram of Figure~\ref{F:Fig2} suffices
to let one disk appear.

\section{Disks and removable pairs of letters}
\label{S:Disk}

\begin{figure} [htb]
$$\includegraphics{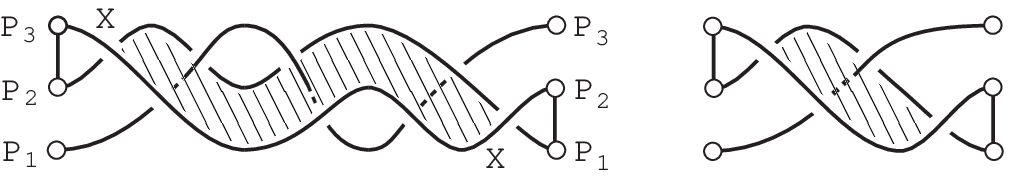}$$
\caption{\smaller A $(2, 1)$-disk (left), and a diagram
that is {\it not} a $(2, 1)$-disk (right): the third strand
pierces the ribbon made by the first two strands;
it is convenient in the formal definition to
appeal to the points~$P_i$ and~$P'_j$, but, in
essence, the disk is the part lying between the
crossings denoted~$X$ and~$X'$.}
\label{F:Fig1}
\end{figure}

\begin{defi}
(Figure~\ref{F:Fig1}) Assume that $\diag$ is an
$n$-strand braid diagram, which is the projection of
a 3-dimensional geometric braid~$\b$ consisting
of $n$~disjoint curves connecting $n$~points
$P_1, \pp P_n$ in the plane $z = 0$ to
$n$~points $P'_1, \pp, P'_n$ in the plane $z = 1$.
For $1 \le i, j < n$, we say that $\diag$ is an {\it $(i,
j)$-disk} if $\diag$ begins with a crossing of the
strands starting at~$P_i$ and~$P_{i+1}$, it finishes
with a crossing of opposite orientation of the
strands ending at~$P'_j$ and~$P'_{j+1}$, and the
figure obtained from~$\b$ by connecting
$P_i$ to~$P_{i+1}$ and $P'_j$ to~$P'_{j+1}$ is
isotopic to the union of $n-2$~curves and the
boundary of a disk disjoint from these curves.
\end{defi}

This definition is directly reminiscent of the notion of
a {\it life disk} in a singular braid introduced
in~\cite{FKR}: another way to state that $\diag$ is a
disk is to say that, when one makes the initial and
the final crossings in~$\diag$ singular---with the
convention that the first crossing is replaced with a
``birth'' singular crossing, while the last one, which is
supposed to have the opposite orientation, is replaced
with a ``death'' singular crossing---then the resulting
figure is a life disk.

We shall address Question~\ref{Q:1} using the braid
group~$B_n$ and the geometry of its Cayley graph. As
is standard, braid diagrams will be encoded by finite
words over the alphabet
$\{\s_1^{\pm 1}, \pp, \s_{n-1}^{\pm 1}\}$, using
$\s_i$ to encode the elementary diagram where the
$(i+1)$th strand crosses over the $i$th strand. For
instance, the first three diagrams above are coded by
$\ss1 \sss1$,
$\ss1 \ss2 \ss1 \sss2 \sss1 \sss2$, and $\ss1
\s_2^2\ss1\s_2^2\s_1^{-2}\sss2\s_1^{-2}\sss2$,
respectively. 

We denote by~$\equiv$ the equivalence relation on
braid words that corresponds to braid isotopy. As is well
known, $\equiv$ is the congruence generated by the
pairs $(\ss i \ss j, \ss j \ss i)$ with $\vert i - j\vert \ge
2$ and $(\ss i \ss j \ss i, \ss j
\ss i \ss j)$ with $\vert i - j\vert = 1$, together with
$(\ss i\sss i, \e)$ and $(\sss i \ss i, \e)$, where
$\e$ denotes the empty word.

\begin{prop}
A braid diagram is an $(i, j)$-disk if and only if it is
encoded in a word of the form~$\s_i^e w \s_j^{-e}$
with $e = \pm 1$ and $\s_i^e w \s_j^{-e} \equiv w$.
\end{prop}

\begin{proof}
Assume that $\diag$ is an $(i, j)$-disk. By definition,
$\diag$ is encoded in some braid word of the form
$\s_i^e w \s_j^{-e}$ with $e = \pm 1$. Moreover, we
can assume that, after an isotopy, the strands of~$D$
starting at positions~$i$ and~$i+1$ make an unpierced
ribbon. Then, the initial $\s_i^e$ crossing may be
pushed along that ribbon, so as to eventually cancel the
final $\s_j^{-e}$ crossing. Hence $D$ is
isotopic to the diagram obtained by deleting its first
and last crossings,, \ie, we have $\s_i^e
w \s_j^{-e} \equiv w$.

Conversely, assume that $D$ is encoded in $\s_i^e w
\s_j^{-e}$ and $\s_i^e w \s_j^{-e} \equiv w$ holds.
Then we have $\s_i^e w \equiv w \s_j^e$. By
Theorem~2.2 of~\cite{FRZ}, this implies that $D$
contains a ribbon connecting $[i, i+1] \times 0$ to $[j,
j+1] \times 1$ that is, up to an isotopy, disjoint from the
other strands. Hence, with our current definition, $D$ is
an $(i, j)$-disk.
\end{proof}

We thus are led to introduce:

\begin{defi}
A braid word of the form $\s_i^e w \s_j^{-e}$ with $e =
\pm 1$ is said to be a {\it removable pair of letters} if
$\s_i^e w \s_j^{-e} \equiv w$ holds. 
\end{defi}

With this notion, Question~\ref{Q:1} is equivalent to

\begin{ques}\label{Q:2}
Does every nonempty trivial braid word contain a
removable pair of letters?
\end{ques}

Speaking of ``removable pair'' is natural here: indeed,
saying that a braid word~$w'$ contains a removable pair
$\s_i^e w \s_j^{-e}$ implies that $w'$ is equivalent to
the word obtained from~$w'$ by replacing the subword
$\s_i^e w \s_j^{-e}$ with~$w$, \ie, by deleting the end
letters~$\s_i^e$ and~$\s_j^{-e}$. Observe that the
notion of a removable pair of letters actually makes
sense for any group presentation: we shall use it in a
more general context in Section~\ref{S:Pairs} below.

\begin{figure} [htb]
$$\includegraphics[scale=0.7]{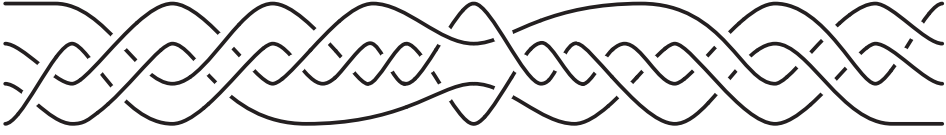}$$
\caption{\smaller A trivial 4-strand braid
diagram containing no disk}
\label{F:Fig2}
\end{figure}

As there exist efficient algorithms for deciding braid
word equivalence, it is easy to systematically search
the possible removable pairs in  a braid word, and an
experimental approach of Question~\ref{Q:2} is
possible. Random tries would suggest a positive answer,
but this is misleading: for instance, the
4~strand braid word
$$\sss1 \s_2^{-2} \sss3 \sss1 \s_2^{-2} \sss1 \sss3
\s_2^{-2} \sss3 \s_2^{-3} \s_1^{-2}
\s_3^2\s_2^3 \ss1 \s_2^2
\ss1\ss3 \s_2^2 \ss3\ss1 \s_2^2 \ss3$$
contains no removable pair of letters, and,
therefore, the associated braid diagram, which is
diplayed in Figure~\ref{F:Fig2}, contains no
disk. This establishes the negative part of
Proposition~\ref{P:Pair}.

\section{The valuation of a pure simple element}

The proof of (the positive part of)
Proposition~\ref{P:Pair} relies on partitioning
the Cayley graph of~$B_3$ using integer
parameters connected with division in the braid
monoid~$B_3^+$. The construction is not specific
to the braid group~$B_3$, nor is it either
specific to braid groups: actually, it is relevant for
all spherical type Artin--Tits groups, and, more
generally, for all Garside groups in the sense
of~\cite{Dgk}.

A monoid~$\GGp$ is said to be a {\it Garside
monoid} if it is cancellative, $1$ is the
only invertible element, any two elements admit a
left and a right least common multiple, and $\GGp$
contains a Garside element, defined as an element
whose left and right divisors coincide, they
generate the monoid, and they are finite in
number. If $\GGp$ is a Garside monoid, it embeds in
a group of fractions. A group~$G$ is said to be a
{\it Garside group} if
$G$ can be expressed in at least one way as the
group of fractions of a Garside monoid. 

Typical examples of Garside monoids are the braid
monoids~$B_n^+$, and, more generally, the Artin--Tits 
monoids~$A^+$ of spherical type, \ie, those Artin--Tits
monoids such that the associated Coxeter group~$W$ is
finite. In this case, the image of the longest element
of~$W$ under the canonical section of the projection
of~$A^+$ onto~$W$ is a Garside element in~$A^+$. In
the particular case of~$B_n^+$, one obtains the
half-twist braid~$\D_n$. So, the braid
groups~$B_n$, and, more generally, the Artin--Tits
groups of  sperical type, are Garside groups. Let us
mention that a given group may be the group of
fractions of several Garside monoids: for instance, the
braid groups~$B_n$ admit a second Garside structure,
associated with the Birman--Ko--Lee monoid
of~\cite{BKL}---see \cite{Bes, PicNote} for similar results
involving other Artin--Tits groups. Still another Garside
structure for~$B_3$ involves the submonoid generated
by~$\ss1$ and~$\ss1\ss2$, a Garside monoid with
presentation $\langle a, b ; aba = b^2\rangle$, hence
not of Artin--Tits type.

Assume that $\GGp$ is a Garside monoid. Then
every element~$x$ in~$\GGp$ admits finitely many
expressions as a product of atoms (indecomposable
elements), and the supremum~$\Vert x\Vert$ of the
length of these decompositions, called the {\it
norm} of~$x$, satisfies $\Vert xy\Vert \ge \Vert
x\Vert + \Vert y\Vert$ and $\Vert x\Vert \ge 1$ for
$x \not= 1$. Then there exists in~$\GGp$ a unique
Garside element of minimal norm; this element is
traditionally denoted~$\D$, and its (left and right)
divisors are called the {\it simple} elements
of~$\GGp$.

We shall start from two technical results about division
in Garside monoids---as shown in~\cite{Dgk}, these
results also happen to be crucial in the construction of
an automatic structure \cite{Eps, Cha, Chb}. For
$\xx, \yy$ in a Garside monoid~$\GGp$, we denote by
$\xx \dr \yy$ the unique element~$\zz$ such that $\xx
\zz$ is the right lcm of~$\xx$ and~$\yy$, and we write
$\yy \divel \zz$ (\resp $\zz \multer \yy$) to express
that $\yy$ is a left (\resp right) divisor of~$\zz$.

\begin{lemm}\label{L:cover}
Assume that $\GGp$ is a Garside monoid, that $\yy,
\zz$ are elements of~$\GGp$ and that every simple
right divisor of~$\yy\zz$ is a right divisor of~$\zz$. Let
$\xx$ be an arbitrary element of~$\GGp$, and let
$\yy' = \xx \dr \yy$ and $\zz' = (\yy \dr \xx) \dr
\zz$. Then every simple right divisor of~$\yy'\zz'$
is a right divisor of~$\zz'$. 
\end{lemm}

\begin{proof}
Let $\xx' = \yy \dr \xx$ and $\xx'' = \zz \dr (\yy \dr
\xx)$. By definition of a right lcm, we have $\xx \yy'
= \yy \xx'$, and $\xx' \zz' = \zz \xx''$. Moreover
$1$ is the only common right divisor of~$\yy'$
and~$\xx'$. Assume that $\sp$ is a simple right
divisor of~$\yy' \zz'$. Then we have $\xx \yy' \zz'
\multer \sp$, hence $\yy \zz \xx'' \multer \sp$. Let
$\sp' \xx''$ be the left lcm of~$\sp$ and~$\xx''$. Then
$\yy \zz \xx'' \multer \sp$ implies $\yy \zz \xx''
\multer \sp' \xx''$, hence $\yy \zz \multer \sp'$.
Moreover, $\sp$ being simple implies that
$\sp'$ is simple as well, as shows an induction on
the minimal number~$p$ such that $\xx''$ can be
decomposed into the product of $p$~simple
elements. Then, the hypothesis of the lemma
implies $\zz \multer \sp'$, and, therefore, $\zz
\xx'' \multer \sp$, \ie, $\xx' \zz' \multer \sp$. It follows
that $\sp$ is a right divisor of the right lcm of~$\yy'
\zz'$ and~$\xx'
\zz'$, which is $\zz'$ since $1$ is the only
common right divisor of~$\yy'$ and~$\xx'$. 
\end{proof}

\begin{lemm}\label{P:main}
Assume that $\GGp$ is a Garside monoid, that $\yy,
\zz, \xx$ are elements of~$\GGp$, and that every
simple right divisor of~$\yy\zz$ is a right divisor
of~$\zz$.  Then $\yy \not\divel \xx$ implies $\yy\zz
\not\divel \xx\spp$ for every simple element~$\spp$
of~$\GGp$.
\end{lemm}

\begin{proof}
We assume $\yy\zz \divel \xx\spp$, and aim at
proving $\yy \divel \xx$. Let $\yy' = \xx \dr \yy$, and
$\zz' = (\yy \dr \xx) \dr \zz$. By construction, we have
$\yy' \zz' = \xx \dr (\yy \zz)$, and $\yy\zz \divel \xx
\spp$ implies $\yy' \zz' \divel \spp$, so, in particular,
$\yy' \zz'$ must be simple. By Lemma~\ref{L:cover},
every simple right divisor of~$\yy' \zz'$ is a right
divisor of~$\zz'$, so we deduce $\zz' \multer \yy' \zz'$,
which is possible for $\yy' = 1$ only, \ie, for $\yy \divel
\xx$.
\end{proof}

Now, the idea is to consider, for each element of a
Garside group~$\GG$ and each simple element~$\sp$
of~$\GGp$, the maximal power of~$\sp$ that
divides a given element. We begin with the monoid.

\begin{defi}
Assume that $\GGp$ is a Garside monoid. We say that a
simple element~$\sp$ of~$\GGp$ is {\it pure} if $\sp$
is the maximal simple right divisor of~$\sp^k$, for
every~$k$. If $\sp$ is a pure simple element
of~$\GGp$, we define the {\it (left)
valuation}~$\val_\sp(\xx)$ of~$\sp$ in~$\xx$ to be the
maximal~$k$ satisfying $\sp^k \divel \xx$.
\end{defi}

In the braid monoid~$B_n^+$, each
generator~$\ss i$, as well as the Garside
element~$\D_n$---and, more generally, each simple
braid which is an lcm of generators~$\ss i$---is a pure
simple element. If $\GGp$ is an arbitrary Garside
monoid, the Garside element~$\D$ is always pure by
definition, but the atoms or their lcms need not be pure
in general: for instance, in the monoid $\langle a, b \,;\,
aba = b^2\rangle^+$, the atom~$b$ is not pure, as
$b^2$ is simple.

\begin{lemm}\label{L:ineq}
Assume that $\GGp$ is a Garside monoid and that
$\sp$ is a pure simple element of~$\GGp$. Then, for
every~$\xx$ in~$\GGp$, we have
\begin{equation}\label{E:ineq}
\val_\sp(\xx) \le \val_\sp(\xx \spp) 
\le \val_\sp(\xx) + 1
\end{equation}
whenever $\spp$ is a simple element of~$\GGp$; more
specifically, for $\spp = \D$, we have
\begin{equation}\label{E:equal}
\val_\sp(\xx \D) = \val_\sp(\xx) + 1.
\end{equation}
\end{lemm}

\begin{proof}
First $\sp^k \divel \xx$ implies $\sp^k
\divel \xx \spp$ for every~$\spp$, hence $\val_\sp(\xx) 
\le \val_\sp(\xx \spp)$. On the other hand,  assume
$\sp^{k+1} \not\divel \xx$.   By hypothesis, every
right divisor of~$\sp^{k+2}$ is a right divisor
of~$\sp$. Applying Lemma~\ref{P:main} with
$\yy = \sp^{k+1}$ and $\zz = \sp$, we deduce
$\sp^{k+2} \not\divel \xx \spp$, hence $\val_\sp(\xx
\spp)  \le \val_\sp(\xx) + 1$, and~\eqref{E:ineq}
follows.

As $\D$ is simple, \eqref{E:ineq} implies $\val_\sp(\xx
\D) \le \val_\sp(\xx) + 1$. On the other hand,
let~$\ff$ be the automorphism  of~$\GGp$ defined
for $\zz$ a simple element by $\ff(\zz) = (\zz
\dr \D) \dr \D$ (see~\cite{Dgk}). Then $\zz \D = \D
\ff(\zz)$ holds for every~$\zz$. Now assume $\sp^k
\divel \xx$. We find
$$\xx \D = \sp^k \xx' \D = \sp^k \D \ff(\xx') =
\sp^{k+1} (\sp \dr \D) \ff(\xx'),$$
hence $\sp^{k+1} \divel \xx \D$, and, therefore, 
$\val_\sp(\xx \D) > \val_\sp(\xx)$,
hence~\eqref{E:equal}.
\end{proof}

We now extend the maps~$\val_\sp$
from a Garside monoid~$\GGp$ to its group of
fractions~$\GG$. As $\D$ is a common multiple of all
atoms in~$\GGp$, every element of~$\GG$ can
be expressed as~$\xx \D^k$ with~$\xx \in
\GGp$ and $k \in \ZZ$. Unless we require that $k$ be
maximal, the decomposition need not be unique.
However, we have the following result:

\begin{lemm}
Assume that $\GGp$ is a Garside monoid, $\xx,
\xx'$ are elements of~$\GGp$, and we have
$\xx \D^k = \xx' \D^{k'}$ in the group of
fractions~$\GG$ of~$\GGp$. Then, for each pure
simple element~$\sp$ of~$\GGp$, we have
$\val_\sp(\xx) + k = \val_\sp(\xx') + k'$.
\end{lemm}

\begin{proof}
Assume for instance $k \le k'$, say $k' = k + m$.
Then we have $\xx \D^k = \xx' \D^m \D^k$ in~$\GG$,
hence $\xx = \xx' \D^m$ in~$\GGp$ (we recall that
$\GGp$ embeds in~$\GG$). Using Lemma~\ref{L:ineq}
$m$~times, we obtain
$\val_\sp(\xx'
\D^m) = \val_\sp(\xx') + m$ for every~$\sp$, hence
$\val_\sp(\xx) = \val_\sp(\xx') + m$, \ie, 
$\val_\sp(\xx) + k = \val_\sp(\xx') + k'$.
\end{proof}

Then the following definition is natural:

\begin{defi}
Assume that $\GGp$ is a Garside monoid, $\GG$ is the
group of fractions of~$\GGp$, and $\sp$ is a pure
simple element of~$\GGp$. Then, for~$\xx$ in~$\GG$,
the {\it (left) valuation}~$\val_\sp(\xx)$ of~$\sp$
in~$\xx$ is defined to be $\val_\sp(\zz) + k$, where $\xx
= \zz \D^k$ is an arbitrary decomposition of~$\xx$ with
$\zz \in \GGp$ and $k \in \ZZ$.
\end{defi}

\begin{exam}\label{X:exp}
Let $G = B_3$, and $\xx = \sss1 \ss2$. We can also write
$\xx = \ss2 \s_1^2 \D_3^{-1}$. We have
$\val_{\ss1}(\ss2 \s_1^2) = 0$ and $\val_{\ss2}(\ss2
\s_1^2) = 1$, so we find $\val_{\ss1}(\xx) = 0 -1 = -1$,
and 
$\val_{\ss2}(\xx) = 1 -1 = 0$. 
\end{exam}

It is now easy to see that the inequalities of
Lemma~\ref{L:ineq} remain valid in the group:

\begin{prop}\label{P:exp}
Assume that $\GG$ is the Garside group associated with a
Garside monoid~$\GGp$, and that  $\sp$ is a pure
simple element of~$\GGp$. Then, for
every element~$\xx$ in~$\GG$, and every simple
element~$\spp$ in~$\GGp$, we have
\begin{equation}\label{E:ineqg}
\val_\sp(\xx) \le \val_\sp(\xx \spp) 
\le \val_\sp(\xx) + 1;
\end{equation}
for $\spp = \D$, we have $\val_\sp(\xx\D) =
\val_\sp(\xx) + 1$.
\end{prop}

\begin{proof}
Assume $\xx = \yy \D^k$ with $\yy \in \GGp$. We
have $\xx \spp = \yy \D^k \spp = \yy \ff^{-k}(\spp)
\D^k$. Then $\yy \ff^{-k}(t)$ belongs
to~$\GGp$, hence we have 
$\val_\sp(\xx) = \val_\sp(\yy) + k$, and
$\val_\sp(\xx \spp) = \val_\sp(\yy \ff^{-k}(\spp)) + k$.
As $\ff^{-k}(\spp)$ is a simple element
of~$\GGp$, Lemma~\ref{L:ineq} gives
$$\val_\sp(\yy) \le \val_\sp(\yy \ff^{-k}(\spp))
\le \val_\sp(\yy) + 1,$$
so \eqref{E:ineqg} follows. The result  for $\spp = \D$ is
obvious, since we obtain $\xx \D = \yy
\D^{k+1}$, hence $\val_\sp(\xx \D) = \val_\sp(\yy) +
k +1 = \val_\sp(\xx) + 1$ directly.
\end{proof}

Inequality~\eqref{E:ineqg} is the algebraic socle on
which we shall build in the sequel.

\section{Partitions of the Cayley graph}\label{S:Part}

From now on, we restrict to  Artin--Tits groups, \ie, we
consider presentations of the form
\begin{equation}\label{E:Pres}
\langle S \,;\, \Prod(\gen, \genn, m_{\gen,\genn}) =
\Prod(\genn, \gen, m_{\gen, \genn}) \mbox{ for $\gen \not = \genn$
in~$S$ }\rangle, 
\end{equation}
where $\Prod(\gen, \genn, m)$ denotes the alternated product
$\gen \genn \gen \genn\ldots$ with $m$~factors, and $m_{\gen, \genn}
\ge 2$ holds. Moreover, we restrict to the spherical
type, \ie, we assume that the Coxeter group obtained by
adding to~\eqref{E:Pres} the relation $\gen^2 = 1$ for
each~$\gen$ in~$S$ is finite. Then the monoid~$\AAp$
defined by~\eqref{E:Pres} is a Garside monoid, and the
group~$\AA$ defined by~\eqref{E:Pres} is the group of
fractions of~$\AAp$.

In this case, each generator~$\gen$ in~$S$ is pure, since
$\gen^2$ is not simple and $\gen$ is the right gcd of~$\gen^2$
and~$\D$. Hence, each element~$\xx$ of the
group~$\AA$ has a well-defined
valuation~$\val_\gen(\xx)$ for each~$\gen$ in~$S$,
and we can associate to~$\xx$ the valuation
sequence~$(\val_\gen(\xx); \gen \in \SS)$.

\begin{exam}
Consider the case of ~$B_3$. There are two atoms,
namely~$\ss1$ and~$\ss2$. The valuation sequence
associated with~$\ss1$ is $(1, 0)$, while the one
associated with $\sss1\ss2$ is~$(1, -1)$, as was seen
above. Observe that the influence of right multiplication
on the valuation sequence may be anything that is
compatible with the constraints of~\eqref{E:ineqg}. For
instance, $\ss1$, $\ss1 \ss2$, and $\ss1 \s_2^2 \ss1$
all admit the valuation sequence $(1, 0)$, while the
valuation sequences of~$\ss1 \cdot \ss2$, $\ss1\cdot
\ss1$, $\ss1 \ss2 \cdot \ss1$, and $\ss1 \s_2^2 \ss1
\cdot \ss2$ are $(1, 0)$, $(2, 0)$, $(1, 1)$, and $(2, 1)$,
respectively.
\end{exam}

Using the valuation sequence, we can partition
the group~$\AA$, hence, equivalently, its Cayley
graph, into disjoint regions according to the values of
the valuations. For our current purpose, we shall
consider a coarser partition, namely the one obtained by
taking into account not the values of the valuations, but
their relative positions only. Let us say that two
$n$-tuples of integers
$(k_1, \pp, k_n)$ and $(k'_1, \pp, k'_n)$ are {\it
order-equivalent} if $k_i = k_j$ and $k'_i = k'_j$ (\resp
$k_i < k_j$ and $k'_i < k'_j$) hold for the same
pairs~$(i, j)$. The equivalence class of a tuple $(k_1,
\pp, k_n)$ will be called its {\it order-type}. For
instance, there are 3~order-types of pairs,
corresponding to pairs~$(k_1, k_2)$ with $k_1 <
k_2$, $k_1 = k_2$, and $k_1 > k_2$, respectively.
Similarly, there are 13 order-types of triples, and, in the
general case of $n$-tuples, the number of order-types is
the $n$th ordered Bell number $\sum_{p=1}^n a_p
p^n$ with $a_p = \sum_{q = 0}^{n-p} (-1)^q {p+q
\choose q}$.

\begin{defi}
Assume that $\AA$ is an Artin--Tits group of spherical 
type with presentation~\eqref{E:Pres}. For~$\xx$
in~$\AA$, the {\it type} of~$\xx$ is defined to be the
order-type of the sequence $(\val_\s(\xx) \,;\,
\s \in \SS)$.
\end{defi}

So, there are 3~types of braids in~$B_3$, according
to whether the value of~$\val_{\at}$ is smaller than,
equal to, or bigger than the value of~$\val_{\att}$.
These types will be denoted $[\val_{\at} <
\val_{\att}]$, $[\val_{\at} = \val_{\att}]$, and
$[\val_{\at} > \val_{\att}]$. Thus, saying that a
braid~$\b$ in~$B_3$ is of type~$[\val_{\at} >
\val_{\att}]$ means that there are ``more $\ss1$'s than
$\ss2$'s at the left of~$\b$''. For instance, the type
of~$\ss1$ is~$[\val_{\at} > \val_{\att}]$, while that
of~$\ss2$ and of~$\sss1$ is~$[\val_{\at} <
\val_{\att}]$, and that of~$1$
or~$\D_3^k$ is~$[\val_{\at} = \val_{\att}]$.

Proposition~\ref{P:exp} immediately leads to constraints
on how the type may change under right
multiplication by a simple element:

\begin{prop}\label{P:neighbour}
Assume that $\AA$ is an Artin--Tits group of spherical type. 
Say that two types~$T, T'$ are {\it neighbours}
if there exist $(k_1, \pp, k_n)$ in~$T$ and $(k'_1, \pp,
k'_n)$ in~$T'$ such that $k'_i - k_i$ is either~$0$
or~$1$ for every~$i$, or is either~$0$ or~$-1$ for
every~$i$. Then, for every~$\xx$ in~$\AA$ and every
simple element~$\spp$ of~$\AAp$, the type
of~$\xx\spp^{\pm 1}$ is a neighbour of the type
of~$\xx$.
\end{prop}

\begin{figure} [htb]
$$\includegraphics*[scale=0.9]{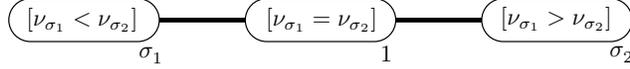}$$
\caption{\smaller The 3 types of braids in $B_3$}
\label{F:Fig3}
\end{figure}

\begin{figure} [htb]
$$\includegraphics*[scale=0.9]{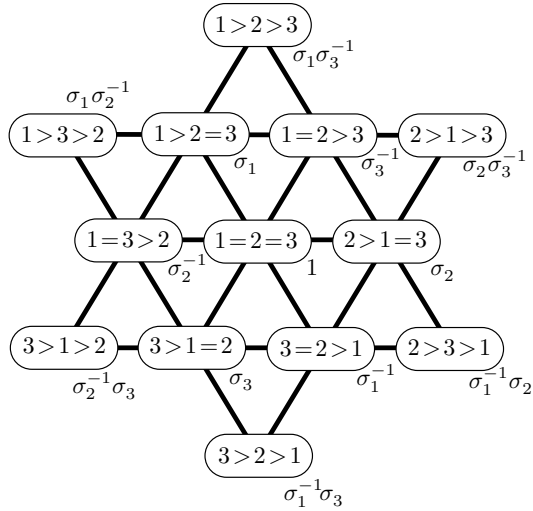}$$
\caption{\smaller The 13 types of braids in
$B_4$---here $i$ stands for $\val_{\s_i}$}
\label{F:Fig4}
\end{figure}

We display in Figures~\ref{F:Fig3} and~\ref{F:Fig4} the
graph of the neighbour relation for order-types of pairs
and of triples---as well as examples of $3$- and
$4$-strand braids of the corresponding types.
We see in Figure~\ref{F:Fig3} that the
types~$[\val_{\at} >
\val_{\att}]$ and~$[\val_{\at} <\val_{\att}]$ are not
neighbours, since, starting with a pair~$(k_1, k_2)$
with~$k_1 > k_2$ and adding~$1$ to~$k_1$ or~$k_2$,
we can obtain $(k'_1, k'_2)$ with $k'_1 \ge k'_2$, but
not with $k'_1 < k'_2$. As a consequence, we cannot
obtain a braid of type~$[\val_{\at} < \val_{\att}]$ by
multiplying a braid of type~$[\val_{\at} > \val_{\att}]$
by a single simple braid or its inverse: crossing the
intermediate type~$[\val_{\at} = \val_{\att}]$ is
necessary. Similarly, we can see on Figure~\ref{F:Fig4}
that, for instance, going from type~$[\val_{\ss1} >
\val_{\ss2} = \val_{\ss3}]$ to type~$[\val_{\ss1} <
\val_{\ss2} = \val_{\ss3}]$ necessitates that one goes
through at least one of the intermediate types
$[\val_{\ss1} = \val_{\ss2} < \val_{\ss3}]$, 
$[\val_{\ss1} = \val_{\ss2} = \val_{\ss3}]$,  or
$[\val_{\ss1} = \val_{\ss3} < \val_{\ss2}]$.

\section{Loops in the Cayley graph}\label{S:Pairs}

We are now ready to establish that every nonempty
trivial $3$-strand braid word contains at least one
removable pair of letters. The geometric idea of the
proof is as follows: a trivial word corresponds to a loop
in the Cayley graph of~$B_3$, and we can choose the
origin of that loop so that it contains vertices of
types~$[\val_{\at} < \val_{\att}]$ and~$[\val_{\at}
> \val_{\att}]$. But then Proposition~\ref{P:neighbour}
tells us that one cannot jump from the
region~$[\val_{\at} < \val_{\att}]$ to the
region~$[\val_{\at} > \val_{\att}]$ without crossing the
separating region, \ie,~$[\val_{\at} = \val_{\att}]$. This
means that some subword of~$w$ must represent a
power of~$\D_3$, and it is easy to deduce a removable
pair of letters.

Actually, we shall prove a more general statement valid
for every Artin--Tits group with two generators, \ie, for
every Artin--Tits group of type~$I_2(m)$---the case
of~$B_3$ corresponding to $m = 3$:

\begin{prop}\label{P:Main}
Assume that $\AA$ is an Artin--Tits group of
type~$I_2(m)$, \ie, $\AA$ admits the presentation
$\langle \at, \att \,;\, \at\att\at\att\pp = \att\at\att\at\pp \rangle$
where both sides of the equality have length~$m$.
Then every nonempty word on the letters~$\at^{\pm
1}, \att^{\pm 1}$ representing~$1$ in~$\AA$ contains a
removable pair of letters.
\end{prop}

We begin with two auxiliary results. 

\begin{lemm}\label{L:Conj}
Assume that $\GG$ is a group generated by a set~$\SS$,
that $w$ is a trivial word on~$\SS \cup \SS\ii$ (\ie, $w$
represents~$1$ in~$\GG$), and some cyclic conjugate
of~$w$ contains a removable pair of letters. Then $w$
contains a removable pair of letters.
\end{lemm}

\begin{proof}
Assume that we have $w = uv$  and $\gen^e w'
\genn^{-e}$ is a removable pair of letters in~$vu$, with
$\gen, \genn \in \SS$, and $e = \pm 1$. Let us write $vu
= w_1 \gen^e w' \genn^{-e} w_2$. If $w_1 \gen^e w'
\genn^{-e}$ is a prefix of~$v$, or if $\gen^e w'
\genn^{-e} w_2$ is a suffix of~$u$, then $\gen^e w'
\genn^{-e}$ is a subword of~$w$, and the result is
obvious. Otherwise, we have $w' = v' u'$ with $\gen^e
v'$ a suffix of~$v$ and $u' \genn^{-e}$ a prefix of~$u$,
hence $v = w_1 \gen^e v'$ and $u = u' \genn^{-e} w_2$.
By construction, $\genn^{-e} w_2 w_1 \gen^e$ is a
subword of~$uv$, \ie, of~$w$. Let us use $\equiv$ for
the congruence that defines~$G$. By hypothesis, we
have $uv \equiv \e$ and  $\gen^e v'u' \genn^{-e}
\equiv v' u'$, hence $\genn^e {u'}\ii {v'}\ii \gen^{-e}
\equiv {u'}\ii {v'}\ii$. We deduce
\begin{align*}
w_2 w_1
\equiv \genn^e {u'}\ii u v {v'}\ii \gen^{-e}
&\equiv \genn^e {u'}\ii {v'}\ii \gen^{-e}\\
&\equiv {u'}\ii {v'}\ii 
\equiv {u'}\ii u v {v'}\ii 
\equiv \genn^{-e} w_2 w_1 \gen^e,
\end{align*}
which shows that $\genn^{-e} w_2 w_1 \gen^e$ is a
removable pair of letters in~$w$.
\end{proof}

\begin{lemm}\label{L:Delta}
Let $\AA$ be an Artin--Tits group with
presentation~\eqref{E:Pres}. Assume that
$\gen, \genn$ belong to~$\SS$, $\cc$ belongs
to~$\SS \cup \SS\ii$, and $w$ is a word on $\SS \cup
\SS\ii$ such that $\genn w \cc \equiv \Prod(\gen,
\genn, m_{\gen, \genn})^k$ holds and $\genn w$
represents an element of the region $[\val_\gen <
\val_\genn]$. Then either $\gen\ii \genn w \cc$ or
$\genn w \cc$ is a removable pair of letters.
\end{lemm}

\begin{proof}
For $u$ a word on~$\SS \cup \SS\ii$, let $\cl u$ denote
the element of~$\AA$ represented by~$u$. Let us write
$m$ for~$m_{\gen, \genn}$. By hypothesis, we have
$\val_\gen(\cl{\genn w \cc}) = \val_\genn(\cl{\genn w
\cc}) = k$. Assume first that $mk$ is even. Then there
are two possibilities for~$\cc$ only, namely $\cc =
\gen$, and $\cc = \genn\ii$. Indeed, $\cl{\genn w
\cc}$ is $\Prod(\gen, \genn, m)^k$, so
$\cc = \gennn^{\pm 1}$ with $\gennn
\not= \gen, \genn$ would imply
$$\val_\gen(\cl{\genn w }) =
\val_\gen(\cl{\genn w \cc}) =
\val_\genn(\cl{\genn w \cc}) =
\val_\genn(\cl{\genn w}),$$
while $\cc = \gen\ii$ and $\cc = \genn$ would imply
$$\val_\gen(\cl{\genn w }) \ge
\val_\gen(\cl{\genn w \cc}) =
\val_\genn(\cl{\genn w \cc}) \ge
\val_\genn(\cl{\genn w}),$$
all contradicting the hypothesis $\val_\gen(\cl{\genn w }) <
\val_\genn(\cl{\genn w})$.

Now, for $\cc = \gen$, we find
$$\gen\ii \genn w  \gen
\equiv \gen\ii \, \Prod(\gen, \genn, m)^k 
\equiv
\Prod(\gen, \genn, m)^k \, \gen\ii 
\equiv \genn w \gen \gen\ii 
\equiv \genn w,$$
\ie, $\gen\ii \genn w \cc$ is a removable pair. Similarly,
for $\cc = \genn\ii$, we find
$$\genn w \genn\ii 
\equiv \Prod(\gen, \genn, m)^k 
\equiv \genn\ii \Prod(\gen, \genn, m)^k \genn
\equiv \genn\ii \genn w \genn\ii \genn
\equiv w,$$
\ie, $\genn w \cc$ is a removable pair. The argument is
similar when $mk$ is odd, the possible values of~$\cc$
now being $\gen\ii$ and~$\genn$ instead of~$\gen$
and~$\genn\ii$.
\end{proof}

\begin{figure}[htb]
$$\includegraphics{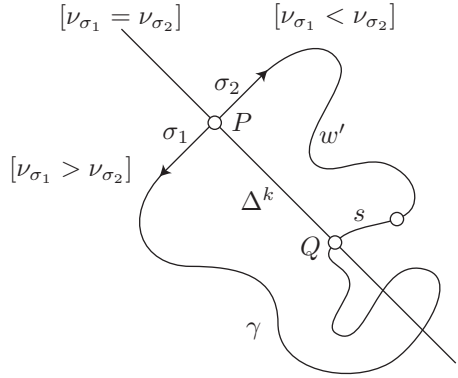}$$
\caption{Proof of Proposition~\ref{P:Main}: a loop must
intersect the diagonal at least twice}
\label{F:Fig5}
\end{figure}

\begin{proof}[Proof of Proposition~\ref{P:Main}]
(Figure~\ref{F:Fig5})
Assume that $w$ is a nonempty word on the letters
$\at^{\pm1}$, $\att^{\pm1}$ representing~$1$.
Necessarily $w$ contains the same number of letters
with exponent~$+1$ and with exponent~$-1$, so it must
contain a subword of the form~$\cc\ii \ccc$ or $\cc
\ccc\ii$ with $\cc, \ccc \in \{\at, \att\}$. Assume for
instance that $w$ contains a subword of the
form~$\cc\ii  \ccc$; the argument in the case of~$\cc
\ccc\ii$ would be similar. The case $\cc = \ccc$ is trivial
(then $\cc\ii \cc$ is a removable pair of letters of~$w$,
and we are done), so, up to a symmetry, we can assume
that $\cc\ii \ccc$ is~$\at\ii \att$. 

The word~$w$ specifies a path~$\g$ in the Cayley graph
of~$G$, and, by hypothesis, $\g$ is a loop. Let $P$ be
the point of~$\g$ corresponding to the middle vertex in
the subword~$\at\ii \att$ considered above. By
Lemma~\ref{L:Conj}, we can assume that $P$ is the origin 
of~$\g$ without loss of generality.

Now, let us follow~$\g$ starting from~$P$: as the first
letter is~$\att$, the path~$\g$ enters the
region~$[\val_{\at} < \val_{\att}]$. At the other end, the
last letter of~$\g$ is~$\at\ii$, which means that, before
ending at~$P$, the path~$\g$ comes from the
region~$[\val_{\at} > \val_{\att}]$. So
$\g$ goes from the region~$[\val_{\at} < \val_{\att}]$
to the region~$[\val_{\at} > \val_{\att}]$. By
Proposition~\ref{P:neighbour}, $\g$ must cross the
separating region~$[\val_{\at} = \val_{\att}]$ at least
once. This means that there must exist at least one
second point~$Q$ in~$\g$ with type~$[\val_{\at} =
\val_{\att}]$. Now---and this is where we use the
hypothesis that $\AA$ is of Coxeter type~$I_2(m)$---the
only elements of~$\AA$ of this type are the
powers of the element~$\D$, \ie, of~$\Prod(\at, \att,
m)$. So we deduce that (a cyclic conjugate
of)~$w$ must contain a subword~$\at\ii
\att w' \cc$ such that $\att w' \cc$ is equivalent
to a power of~$\Prod(\at, \att, m)$ and $\att w'$
represents an element of the region $[\val_{\at} <
\val_{\att}]$. Then Lemma~\ref{L:Delta} implies that
either
$\at\ii \att w' \cc$ or $\att w' \cc$ is a removable pair
of letters.
\end{proof}

\begin{rema}
It is known \cite{Bst, BrW, CMW} that the Cayley
graph of any Garside group is traced on some flag
complex of the form ${\mathcal X} \times {\mathbb
R}$, where the
${\mathbb R}$-component corresponds to powers
of~$\D$. In the case of an Artin--Tits group of
type~$I_2(m)$, the space~${\mathcal X}$ is an
$m$-valent tree. A loop~$\g$ in the Cayley graph
projects onto a loop in the tree, so the
projection necessarily goes twice through the same
vertex, which means that $\g$ contains vertices
that are separated by a power of~$\D$, and we can
deduce the existence of a removable pair of
letters as above.
\end{rema}

\section{Special cases}

As the counter-example of Figure~\ref{F:Fig2}
shows, a trivial $4$-strand braid word need not contain
any removable pair of letters. However, partial positive results exist,
in particular when we consider words of the form~$u\ii
v$, with $u, v$ positive words representing a divisor
of~$\D$.

The following result is an easy consequence of the
classical Exchange Lemma for Coxeter groups (\cite{Bou},
Lemma~IV.1.4.3) rephrased for Artin--Tits monoids. 

\begin{lemm}\label{L:Eddy}
Assume that $\AAp$ is an Artin--Tits monoid of spherical type, $\gen, \genn$ are atoms of~$\AAp$, and we have $\gen
\not\divel \xx$ and $\gen \divel \xx\genn \divel \D$. Then
we have $\xx\genn = \gen \xx$.
\end{lemm}

Indeed, let~$\pi$ denote the bijection of the divisors
of~$\D$ in~$\AAp$ to the corresponding
Coxeter group~$W$ and $\ell$ denote the length
in~$W$. Then $\gen \divel \xx \genn$ implies
$\ell(\pi(\gen \xx \genn)) < \ell(\pi(\xx \genn))$.
Hence the minimal decomposition of
$\pi(\gen \xx \genn)$ is obtained from that of~$\pi(\xx
\genn)$ by removing one generator, which cannot
come from~$\xx $ for, otherwise, we would obtain
$\ell(\pi(\gen \xx)) < \ell(\pi(\xx))$ by
cancelling~$\genn$ and contradict $\gen
\not\divel \xx$. So we must have $\pi(\gen \xx \genn) =
\pi(\xx)$, hence $\pi(\gen \xx) = \pi(\xx \genn)$,
in~$W$, and $\gen \xx  = \xx \genn $ in~$\AAp$.

\begin{prop}\label{P:SimpleDisk}
Assume that $\AAp$ is an Artin--Tits monoid of spherical type, and $w$ is a nonempty trivial word of the form
$u\ii v$ with $u, v$ positive simple words. Then $w$
contains at least one removable pair of letters.
\end{prop}

\begin{proof}
For~$w$ a positive word, let~$\cl w$ denote the
element of~$\AAp$ represented by~$w$. Now, let $\gen$
be the first letter in~$u$. By hypothesis, we have
$\gen \divel \cl v$. Let $v'\genn$ be the shortest prefix of~$v$
such that $\gen \divel \cl{v'\genn}$ is true. Then, by definition,
we have $\gen \not\divel \cl{v'}$, and, as $\cl{v}$ is
supposed to be simple, so is $\cl{v'\genn}$. We can
therefore apply Lemma~\ref{L:Eddy}, and we obtain
$\gen v' \equiv v' \genn$, hence $\gen\ii v' \genn \equiv v'$. Thus
$\gen\ii v' \genn$ is a removable pair of letters in~$w$.
\end{proof}

\begin{figure}[htb]
$$\includegraphics*[scale=.85]{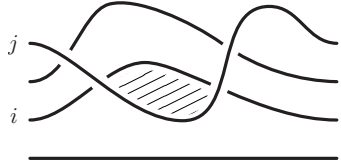}$$
\caption{\smaller Disk in a trivial braid diagram coded
by~$u\ii v$ with $u, v$ simple}
\label{F:Fig6}
\end{figure}

\begin{rema}
In the case of braids, a direct geometric argument
also gives Proposition~\ref{P:SimpleDisk}.
Indeed, if $u$ and $v$ are positive braid words
representing simple braids, then the braid diagrams
coded by~$u$ and~$v$ can be realised as the
projections of three-dimensional figures where the
$i$-th strand entirely lives in the plane~$y = i$: the
simplicity hypothesis guarantees that no altitude
contradiction can occur, as any two strands cross at
most once \cite{Dfq}. So the same is true for the braid
coded by~$u\ii v$, provided we require that the strand
living in the plane
$y = i$ is the one at position~$i$ after~$u\ii$. Now, let
$i$ be the least index such that the $i$th strand is not a
straight line, and let $j$ be the least index such that the
$j$th strand crosses the $i$th strand. Then, necessarily,
the
$i$th and the $j$th strands make a disk, as they must
return to their initial position if $u\ii v$ represents~$1$
(Figure~\ref{F:Fig6}).
\end{rema}

Proposition~\ref{P:SimpleDisk} does not extend to
arbitrary trivial negative--positive words, \ie,  of the
form~$u\ii v$ with $u, v$ positive: the hypothesis that
$u$ and $v$ represent a simple braid is essential. 

An easy method for producing equivalent positive braid
words is as follows: starting with a seed consisting of two
positive words $u, v$, we can complete them into
equivalent words---\ie, we can find a common right
multiple for~$\cl u$ and~$\cl v$---by using the word
reversing technique of~\cite{Dff}, which gives two
positive words~$u', v'$ so that both~$uv'$ and~$vu'$
represent the right lcm of~$\cl u$ and~$\cl v$. Then, by
construction, ${v'}\ii u\ii v u'$ represents~$1$. By
systematically enumerating all possible seeds~$(u, v)$,
we obtain a large number of negative--positive trivial
braid words in which possible removable pair of letters
can be investigated. 

One obtains in this way very few counter-examples, \ie,
trivial braid words with no removable pair of letters. In
the case of~$B_4$, there exists no counter-example with
seeds of length at most~$4$, and there exists only one
counter-example among the
$29,403$ pairs of length~$5$ words, namely the one of
Figure~\ref{F:Fig2}, which is associated with the
seed $(\s_1^2\s_2^3, \s_3^2\s_2^3)$. The situation 
is similar with longer seeds, and for~$B_n$ with $n \ge
5$. This explains why random tries have little chance to
lead to counter-examples, and raises the question of
understanding why there seems to almost always exist
disks in trivial braid diagrams. 

Finally, let us mention a connection with the (open)
question of unbraiding every trivial braid diagram in
such a way that all intermediate diagrams have at most
as many crossings as the initial diagram---as is well
known, there is no solution in the case of knots when
the number of crossings is considered, but there is now
a solution when the complexity is defined in a more
subtle way~\cite{Dyn}. Assume that a method for
detecting removable pairs of letters has been choosen.
Then one obtains an unbraiding algorithm by
starting with an arbitrary braid word and iteratively
removing removable pairs of letters until no one is left.
If the answer to Question~\ref{Q:1} were positive, this
algorithm would always succeed, in the sense that it
would end with the empty word if and only if the initial
word is trivial. Note that the number of iteration steps is
always bounded by half the length of the initial word. In
the case of $4$~strands and more, the answer to
Question~\ref{Q:1} is negative, so the above algorithm
is {\it not} correct. In addition, it must be kept in mind
that, in any case, the algorithm requires a subroutine
detecting removable pairs: we can appeal to any
solution of the braid word problem, but, then, the
algorithm gives no new solution to that word problem,
nor does it either answer the question of
length-decreasing unbraiding as long as there is no
length-preserving method for proving an equivalence of
the form $\s_i^e w \s_j^{-e} \equiv w$.

\begin{figure} [htb]
$$\includegraphics*[scale=.65]{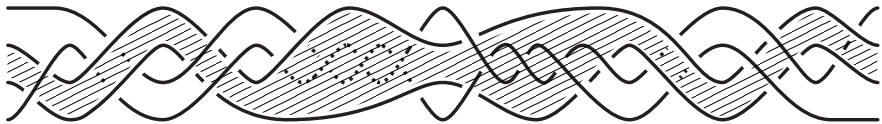}$$
\caption{\smaller A ribbon in the counter-example of
Figure~\ref{F:Fig2}}
\label{F:Fig7}
\end{figure}

As trivial diagrams without disk seem to be
rare, it might happen that, in some sense to be made
precise, the above method almost always works. It can
be observed on Figure~\ref{F:Fig7} that the
braid diagram of Figure~\ref{F:Fig2}, which
contains no disk, contains an actual ribbon, in the
sense that no isotopy is needed to let this ribbon
appear. By merging the two strands bordering this
ribbon, one obtains a $3$-strand diagram---namely
the last example in Section~\ref{S:Intro}---which
contains a disk. Improving the unbraiding method
so as to include such a strand merging procedure
might make it work for still more cases.

\end{document}